\documentclass[a4paper,11pt]{amsart}
\usepackage{amsmath,amssymb,amsfonts,latexsym}

\newtheorem{theorem}{Theorem}[section]

\newtheorem{corollary}[theorem]{Corollary}

\input{tcilatex}

\begin{document}
\title[\textbf{(}$h,q$\textbf{)}-\textbf{Zeta type function with weight}]{%
\textbf{A note on the (}$h,q$\textbf{)}-\textbf{Zeta type function with
weight }$\alpha $}
\author[\textbf{E. Cetin}]{\textbf{Elif Cetin}}
\address{\textbf{Uludag University, Faculty of Arts and Science, Department
of Mathematics, Bursa, Turkey}}
\email{\textbf{elifc2@hotmail.com} }
\author[\textbf{M. Acikgoz}]{\textbf{Mehmet Acikgoz}}
\address{\textbf{University of Gaziantep, Faculty of Science and Arts,
Department of Mathematics, 27310 Gaziantep, TURKEY}}
\email{\textbf{acikgoz@gantep.edu.tr}}
\author[\textbf{I. N. Cangul}]{\textbf{Ismail Naci Cangul}}
\address{\textbf{Uludag University, Faculty of Arts and Science, Department
of Mathematics, Bursa, Turkey}}
\email{\textbf{ncangul@gmail.com}}
\author[\textbf{S. Araci}]{\textbf{Serkan Araci}}
\address{\textbf{University of Gaziantep, Faculty of Science and Arts,
Department of Mathematics, 27310 Gaziantep, TURKEY}}
\email{\textbf{mtsrkn@hotmail.com}}

\begin{abstract}
The objective of this paper is to derive symmetric property of ($h,q$)-Zeta
function with weight $\alpha $. By using this property, we give some
interesting identities for ($h,q$)-Genocchi polynomials with weight $\alpha $%
. As a result, our applications possess a number of interesting property
which we state in this paper.

\vspace{2mm}\noindent \textsc{2010 Mathematics Subject Classification.}
11S80, 11B68.

\vspace{2mm}

\end{abstract}

\maketitle

%%%%%%%%%%%%%%%%%%%%%%%%%%%%%%%%%%%%%%%%%%%%%%%%%%%%%%%%%%%%%%%%%%%

%%%%%%%%%%%%%%%%%%%%%%%%%%%%%%%%%%%%%%%%%%%%%%%%%%%%%%%%%%%%%%%%%%%

%%%%%%%%%%%%%%%%%%%%%%%%%%%%%%%%%%%%%%%%%%%%%%%%%%%%%%%%%%%%%%%%%%%

\section{\textbf{INTRODUCTION}}

%%%%%%%%%%%%%%%%%%%%%%%%%%%%%%%%%%%%%%%%%%%%%%%%%%%%%%%%%%%%%%%%%%%

Recently, T. Kim has developed a new method by using $q$-Volkenborn integral
(or $p$-adic $q$-integral on $%
%TCIMACRO{\U{2124} }%
%BeginExpansion
\mathbb{Z}
%EndExpansion
_{p}$) which has added a weight to $q$-Bernoulli polynomials and
investigated their properties (see \cite{Kim 5}). He also showed that this
polynomials are closely related to weighted $q$-Bernstein polynomials and
derived novel properties of $q$-Bernoulli numbers with weight $\alpha $ by
using symmetric property of weighted $q$-Bernstein polynomials on the $q$%
-Volkenborn integral (for more details, see \cite{Kim 7}). After, Araci $et$ 
$al$. have introduced weighted ($h,q$)-Genocchi polynomials and so defined ($%
h,q$)-Zeta type function with weight by applying Mellin transformation to
generating function of ($h,q$)-Genocchi polynomials with weight $\alpha $
which interpolates for ($h,q$)-Genocchi polynomials with weight $\alpha $ at
negative integers (for details, see \cite{Araci 4}). In this paper, we also
consider ($h$,$q$)-Zeta type function with weight and derive some
interesting properties.

We firstly list some notations as follows:

Imagine that $p$ be a fixed odd prime. Throughout this work $%
%TCIMACRO{\U{2124} }%
%BeginExpansion
\mathbb{Z}
%EndExpansion
,$ $%
%TCIMACRO{\U{2124} }%
%BeginExpansion
\mathbb{Z}
%EndExpansion
_{p},$ $%
%TCIMACRO{\U{211a} }%
%BeginExpansion
\mathbb{Q}
%EndExpansion
_{p}$ and $%
%TCIMACRO{\U{2102} }%
%BeginExpansion
\mathbb{C}
%EndExpansion
_{p}$ will denote by the ring of integers, the field of $p$-adic rational
numbers and the completion of the algebraic closure of $%
%TCIMACRO{\U{211a} }%
%BeginExpansion
\mathbb{Q}
%EndExpansion
_{p},$ respectively. Also we denote $%
%TCIMACRO{\U{2115} }%
%BeginExpansion
\mathbb{N}
%EndExpansion
^{\ast }=%
%TCIMACRO{\U{2115} }%
%BeginExpansion
\mathbb{N}
%EndExpansion
\cup \left\{ 0\right\} $ and $\exp \left( x\right) =e^{x}.$ Let $v_{p}:%
%TCIMACRO{\U{2102} }%
%BeginExpansion
\mathbb{C}
%EndExpansion
_{p}\rightarrow 
%TCIMACRO{\U{211a} }%
%BeginExpansion
\mathbb{Q}
%EndExpansion
\cup \left\{ \infty \right\} $ ($%
%TCIMACRO{\U{211a} }%
%BeginExpansion
\mathbb{Q}
%EndExpansion
$ is the field of rational numbers) denote the $p$-adic valuation of $%
%TCIMACRO{\U{2102} }%
%BeginExpansion
\mathbb{C}
%EndExpansion
_{p}$ normalized so that $v_{p}\left( p\right) =1$. The absolute value on $%
%TCIMACRO{\U{2102} }%
%BeginExpansion
\mathbb{C}
%EndExpansion
_{p}$ will be denoted as $\left\vert \text{ }.\right\vert $, and $\left\vert
x\right\vert _{p}=p^{-v_{p}\left( x\right) }$ for $x\in 
%TCIMACRO{\U{2102} }%
%BeginExpansion
\mathbb{C}
%EndExpansion
_{p}.$ When one speaks of $q$-extensions, $q$ is considered in many ways,
e.g. as an indeterminate, a complex number $q\in 
%TCIMACRO{\U{2102} }%
%BeginExpansion
\mathbb{C}
%EndExpansion
,$ or a $p$-adic number $q\in 
%TCIMACRO{\U{2102} }%
%BeginExpansion
\mathbb{C}
%EndExpansion
_{p},$ If $q\in 
%TCIMACRO{\U{2102} }%
%BeginExpansion
\mathbb{C}
%EndExpansion
$ we assume that $\left\vert q\right\vert <1.$ If $q\in 
%TCIMACRO{\U{2102} }%
%BeginExpansion
\mathbb{C}
%EndExpansion
_{p},$ we assume $\left\vert 1-q\right\vert _{p}<p^{-\frac{1}{p-1}},$ so
that $q^{x}=\exp \left( x\log q\right) $ for $\left\vert x\right\vert
_{p}\leq 1.$ We use the following notation 
\begin{equation}
\left[ x\right] _{q}=\frac{1-q^{x}}{1-q},\text{ \ }\left[ x\right] _{-q}=%
\frac{1-\left( -q\right) ^{x}}{1+q}  \label{equation 1}
\end{equation}

where we want to note that $\lim_{q\rightarrow 1}\left[ x\right] _{q}=x;$
cf. [1-21].

For a fixed positive integer $d$, set%
\begin{eqnarray*}
X &=&X_{d}=\lim_{\overleftarrow{n}}%
%TCIMACRO{\U{2124} }%
%BeginExpansion
\mathbb{Z}
%EndExpansion
/dp^{n}%
%TCIMACRO{\U{2124} }%
%BeginExpansion
\mathbb{Z}
%EndExpansion
, \\
X^{\ast } &=&\underset{\underset{\left( a,p\right) =1}{0<a<dp}}{\cup }a+dp%
%TCIMACRO{\U{2124} }%
%BeginExpansion
\mathbb{Z}
%EndExpansion
_{p}
\end{eqnarray*}

and%
\begin{equation*}
a+dp^{n}%
%TCIMACRO{\U{2124} }%
%BeginExpansion
\mathbb{Z}
%EndExpansion
_{p}=\left\{ x\in X\mid x\equiv a\left( \func{mod}dp^{n}\right) \right\} ,
\end{equation*}

where $a\in 
%TCIMACRO{\U{2124} }%
%BeginExpansion
\mathbb{Z}
%EndExpansion
$ satisfies the condition $0\leq a<dp^{n}$ (see [1-21]).

The following $q$-Haar distribution is defined by T. Kim 
\begin{equation*}
\mu _{q}\left( x+p^{n}%
%TCIMACRO{\U{2124} }%
%BeginExpansion
\mathbb{Z}
%EndExpansion
_{p}\right) =\frac{q^{x}}{\left[ p^{n}\right] _{q}}
\end{equation*}

for any positive $n$ (see \cite{Kim 8}, \cite{Kim 9}).

Let $UD\left( 
%TCIMACRO{\U{2124} }%
%BeginExpansion
\mathbb{Z}
%EndExpansion
_{p}\right) $ be the set of uniformly differentiable function on $%
%TCIMACRO{\U{2124} }%
%BeginExpansion
\mathbb{Z}
%EndExpansion
_{p}$. We say that $f$ is a uniformly differentiable function at a point $%
a\in 
%TCIMACRO{\U{2124} }%
%BeginExpansion
\mathbb{Z}
%EndExpansion
_{p},$ if the difference quotient 
\begin{equation*}
F_{f}\left( x,y\right) =\frac{f\left( x\right) -f\left( y\right) }{x-y}
\end{equation*}

has a limit $f^{%
%TCIMACRO{\U{b4}}%
%BeginExpansion
{\acute{}}%
%EndExpansion
}\left( a\right) $ as $\left( x,y\right) \rightarrow \left( a,a\right) $ and
denote this by $f\in UD\left( 
%TCIMACRO{\U{2124} }%
%BeginExpansion
\mathbb{Z}
%EndExpansion
_{p}\right) .$ In \cite{Kim 8} and \cite{Kim 9}, the $p$-adic $q$-integral
of the function $f\in UD\left( 
%TCIMACRO{\U{2124} }%
%BeginExpansion
\mathbb{Z}
%EndExpansion
_{p}\right) $ is defined by Kim%
\begin{equation}
I_{q}\left( f\right) =\int_{%
%TCIMACRO{\U{2124} }%
%BeginExpansion
\mathbb{Z}
%EndExpansion
_{p}}f\left( \xi \right) d\mu _{q}\left( \xi \right) =\lim_{n\rightarrow
\infty }\sum_{\xi =0}^{p^{n}-1}f\left( \xi \right) \mu _{q}\left( \xi +p^{n}%
%TCIMACRO{\U{2124} }%
%BeginExpansion
\mathbb{Z}
%EndExpansion
_{p}\right)  \label{equation 2}
\end{equation}

The bosonic integral is considered as the bosonic limit $q\rightarrow 1,$ $%
I_{1}\left( f\right) =\lim_{q\rightarrow 1}I_{q}\left( f\right) $.
Similarly, the $p$-adic fermionic integration on $%
%TCIMACRO{\U{2124} }%
%BeginExpansion
\mathbb{Z}
%EndExpansion
_{p}$\ is defined by Kim \cite{Kim 2} as follows:%
\begin{equation}
I_{-q}\left( f\right) =\lim_{q\rightarrow -q}I_{q}\left( f\right) =\int_{%
%TCIMACRO{\U{2124} }%
%BeginExpansion
\mathbb{Z}
%EndExpansion
_{p}}f\left( x\right) d\mu _{-q}\left( x\right)  \label{equation 3}
\end{equation}

By using fermionic $p$-adic $q$-integral on $%
%TCIMACRO{\U{2124} }%
%BeginExpansion
\mathbb{Z}
%EndExpansion
_{p}$, ($h,q$)-Genocchi polynomials are defined by \cite{Araci 4}%
\begin{eqnarray}
\frac{\widetilde{G}_{n+1,q}^{\left( \alpha ,h\right) }\left( x\right) }{n+1}
&=&\int_{%
%TCIMACRO{\U{2124} }%
%BeginExpansion
\mathbb{Z}
%EndExpansion
_{p}}q^{\left( h-1\right) \xi }\left[ x+\xi \right] _{q^{\alpha }}^{n}d\mu
_{-q}\left( \xi \right)  \label{equation 5} \\
&=&\lim_{n\rightarrow \infty }\frac{1}{\left[ p^{n}\right] _{-q}}\sum_{\xi
=0}^{p^{n}-1}\left( -1\right) ^{\xi }\left[ x+\xi \right] _{q^{\alpha
}}^{n}q^{h\xi }\text{.}  \notag
\end{eqnarray}

\bigskip For $x=0$ in (\ref{equation 5}), we have $\widetilde{G}%
_{n,q}^{\left( \alpha ,h\right) }\left( 0\right) :=\widetilde{G}%
_{n,q}^{\left( \alpha ,h\right) }$ are called ($h,q$)-Genocchi numbers with
weight $\alpha $ which is defined by%
\begin{equation*}
\widetilde{G}_{0,q}^{\left( \alpha ,h\right) }=0\text{ and }q^{h}\frac{%
\widetilde{G}_{m+1}^{\left( \alpha ,h\right) }\left( 1\right) }{m+1}+\frac{%
\widetilde{G}_{m+1}^{\left( \alpha ,h\right) }}{m+1}=\left\{ 
\begin{array}{cc}
\left[ 2\right] _{q}\text{, } & \text{if }m=0, \\ 
0\text{,} & \text{if }m\neq 0\text{.}%
\end{array}%
\right.
\end{equation*}

By (\ref{equation 5}), we have distribution formula for ($h,q$)-Genocchi
polynomials, which is shown by \cite{Araci 4}%
\begin{equation*}
\widetilde{G}_{n+1,q}^{\left( \alpha ,h\right) }\left( x\right) =\frac{\left[
2\right] _{q}}{\left[ 2\right] _{q^{a}}}\left[ a\right] _{q^{\alpha
}}^{n}\sum_{j=0}^{a-1}\left( -1\right) ^{j}q^{jh}\widetilde{G}%
_{n+1,q^{a}}^{\left( \alpha ,h\right) }\left( \frac{x+j}{a}\right) \text{.}
\end{equation*}

By applying some elementary methods, we shall give symmetric properties of
weighted ($h,q$)-Genocchi polynomials and weighted ($h,q$)-Zeta type
function. Consequently, our applications seem to be interesting and
worthwhile for studying in Theory of Analytic Numbers.

\section{\textbf{ON THE (}$h,q$\textbf{)-ZETA-TYPE FUNCTION}}

In this part, we firstly recall the ($h,q$)-Zeta type function with weight $%
\alpha $ which is derived in \cite{Araci 4} as follows:%
\begin{equation}
\widetilde{\zeta }_{q}^{\left( \alpha ,h\right) }\left( s,x\right) =\left[ 2%
\right] _{q}\sum_{m=0}^{\infty }\frac{\left( -1\right) ^{m}q^{mh}}{\left[ m+x%
\right] _{q^{\alpha }}^{s}}  \label{equation 6}
\end{equation}

where $q\in 
%TCIMACRO{\U{2102} }%
%BeginExpansion
\mathbb{C}
%EndExpansion
$, $h\in 
%TCIMACRO{\U{2115} }%
%BeginExpansion
\mathbb{N}
%EndExpansion
$ and $\Re \left( s\right) >1$. It is clear that the special case $h=0$ and $%
q\rightarrow 1$ in (\ref{equation 6}), it reduces to the ordinary
Hurwitz-Euler zeta function. Now, we consider (\ref{equation 6}) in this form%
\begin{equation*}
\widetilde{\zeta }_{q^{a}}^{\left( \alpha ,h\right) }\left( s,bx+\frac{bj}{a}%
\right) =\left[ 2\right] _{q^{a}}\sum_{m=0}^{\infty }\frac{\left( -1\right)
^{m}q^{mah}}{\left[ m+bx+\frac{bj}{a}\right] _{q^{a\alpha }}^{s}}
\end{equation*}

By applying some basic operations to the above identity, that is, for any
positive integers $m$ and $b$, there exist unique non-negative integers $k$
and $i$ such that $m=bk+i$ with $0\leq i\leq b-1$. For $a\equiv 1(\func{mod}%
2)$ and $b\equiv 1(\func{mod}2)$. Thus, we can compute as follows:%
\begin{align}
\widetilde{\zeta }_{q^{a}}^{\left( \alpha ,h\right) }\left( s,bx+\frac{bj}{a}%
\right) & =\left[ a\right] _{q^{\alpha }}^{s}\left[ 2\right]
_{q^{a}}\sum_{m=0}^{\infty }\frac{\left( -1\right) ^{m}q^{mah}}{\left[
ma+abx+bj\right] _{q^{a\alpha }}^{s}}  \label{equation 7} \\
& =\left[ a\right] _{q^{\alpha }}^{s}\left[ 2\right] _{q^{a}}\sum_{m=0}^{%
\infty }\sum_{i=0}^{b-1}\frac{\left( -1\right) ^{i+mb}q^{\left( i+mb\right)
ah}}{\left[ \left( i+mb\right) a+abx+bj\right] _{q^{a\alpha }}^{s}}  \notag
\\
& =\left[ a\right] _{q^{\alpha }}^{s}\left[ 2\right] _{q^{a}}%
\sum_{i=0}^{b-1}\left( -1\right) ^{i}q^{iah}\sum_{m=0}^{\infty }\frac{\left(
-1\right) ^{m}q^{mbah}}{\left[ ab\left( m+x\right) +ai+bj\right] _{q^{\alpha
}}^{s}}  \notag
\end{align}

From this, we can easily discover the following%
\begin{gather}
\sum_{j=0}^{a-1}\left( -1\right) ^{j}q^{jbh}\widetilde{\zeta }%
_{q^{a}}^{\left( \alpha ,h\right) }\left( s,bx+\frac{bj}{a}\right) =
\label{equation 8} \\
\left[ a\right] _{q^{\alpha }}^{s}\left[ 2\right] _{q^{a}}\sum_{j=0}^{a-1}%
\left( -1\right) ^{j}q^{jbh}\sum_{i=0}^{b-1}\left( -1\right)
^{i}q^{iah}\sum_{m=0}^{\infty }\frac{\left( -1\right) ^{m}q^{mbah}}{\left[
ab\left( m+x\right) +ai+bj\right] _{q^{\alpha }}^{s}}  \notag
\end{gather}

Replacing $a$ by $b$ and $j$ by $i$ in (\ref{equation 7}) and so we have the
following%
\begin{equation*}
\widetilde{\zeta }_{q^{b}}^{\left( \alpha ,h\right) }\left( s,ax+\frac{ai}{b}%
\right) =\left[ b\right] _{q^{\alpha }}^{s}\left[ 2\right]
_{q^{b}}\sum_{j=0}^{a-1}\left( -1\right) ^{j}q^{jbh}\sum_{m=0}^{\infty }%
\frac{\left( -1\right) ^{m}q^{mbah}}{\left[ ab\left( m+x\right) +ai+bj\right]
_{q^{\alpha }}^{s}}
\end{equation*}

By considering the above identity in (\ref{equation 8}), we can easily state
the following theorem.

\begin{theorem}
The following%
\begin{equation*}
\frac{\left[ 2\right] _{q^{b}}}{\left[ a\right] _{q^{\alpha }}^{s}}%
\sum_{i=0}^{a-1}\left( -1\right) ^{i}q^{ibh}\widetilde{\zeta }%
_{q^{a}}^{\left( \alpha ,h\right) }\left( s,bx+\frac{bi}{a}\right) =\frac{%
\left[ 2\right] _{q^{a}}}{\left[ b\right] _{q^{\alpha }}^{s}}%
\sum_{i=0}^{b-1}\left( -1\right) ^{i}q^{iah}\widetilde{\zeta }%
_{q^{b}}^{\left( \alpha ,h\right) }\left( s,ax+\frac{ai}{b}\right)
\end{equation*}%
is true.
\end{theorem}

Now, setting $b=1$ in Theorem 2.1, we have the following distribution formula%
\begin{equation}
\widetilde{\zeta }_{q}^{\left( \alpha ,h\right) }\left( s,ax\right) =\frac{%
\left[ 2\right] _{q}}{\left[ 2\right] _{q^{a}}\left[ a\right] _{q^{\alpha
}}^{s}}\sum_{i=0}^{a-1}\left( -1\right) ^{i}q^{ih}\widetilde{\zeta }%
_{q^{a}}^{\left( \alpha ,h\right) }\left( s,x+\frac{i}{a}\right) \text{.}
\label{equation 9}
\end{equation}

If putting $a=2$ in (\ref{equation 9}) leads to the following corollary.

\begin{corollary}
The following identity holds true:%
\begin{equation*}
\widetilde{\zeta }_{q}^{\left( \alpha ,h\right) }\left( s,2x\right) =\frac{%
\left[ 2\right] _{q}}{\left[ 2\right] _{q^{2}}\left[ 2\right] _{q^{\alpha
}}^{s}}\left( \widetilde{\zeta }_{q^{2}}^{\left( \alpha ,h\right) }\left(
s,x\right) -q^{h}\widetilde{\zeta }_{q^{2}}^{\left( \alpha ,h\right) }\left(
s,x+\frac{1}{2}\right) \right) \text{.}
\end{equation*}
\end{corollary}

Taking $s=-m$ into Theorem 2.1, we have the symmetric property of ($h,q$%
)-Genocchi polynomials by the following theorem.

\begin{theorem}
The following identity%
\begin{equation*}
\left[ 2\right] _{q^{b}}\left[ a\right] _{q^{\alpha
}}^{m-1}\sum_{j=0}^{a-1}\left( -1\right) ^{i}q^{ibh}\widetilde{G}%
_{m,q^{a}}^{\left( \alpha ,h\right) }\left( bx+\frac{bi}{a}\right) =\left[ 2%
\right] _{q^{a}}\left[ b\right] _{q^{\alpha }}^{m-1}\sum_{i=0}^{b-1}\left(
-1\right) ^{i}q^{iah}\widetilde{G}_{m,q^{b}}^{\left( \alpha ,h\right)
}\left( ax+\frac{ai}{b}\right)
\end{equation*}%
is true.
\end{theorem}

Now also, setting $b=1$ and replacing $x$ by $\frac{x}{a}$ on the above
theorem, we can rewrite the following ($h,q$)-Genocchi polynomials with
weight $\alpha $.%
\begin{equation*}
\widetilde{G}_{n,q}^{\left( \alpha ,h\right) }\left( x\right) =\frac{\left[ 2%
\right] _{q}}{\left[ 2\right] _{q^{a}}}\left[ a\right] _{q^{\alpha
}}^{n-1}\sum_{i=0}^{a-1}\left( -1\right) ^{i}q^{ih}\widetilde{G}%
_{n,q^{a}}^{\left( \alpha ,h\right) }\left( \frac{x+i}{a}\right) \text{ }%
\left( 2\nmid a\right) \text{.}
\end{equation*}

Due to Araci $et$ $al$. \cite{Araci 4}, we develop as follows 
\begin{eqnarray*}
\sum_{n=0}^{\infty }\widetilde{G}_{n,q}^{\left( \alpha ,h\right) }\left(
x+y\right) \frac{t^{n}}{n!} &=&\left[ 2\right] _{q}t\sum_{m=0}^{\infty
}\left( -1\right) ^{m}q^{mh}e^{t\left[ x+y+m\right] _{q^{\alpha }}} \\
&=&\left[ 2\right] _{q}t\sum_{m=0}^{\infty }\left( -1\right) ^{m}q^{mh}e^{t%
\left[ y\right] _{q^{\alpha }}}e^{\left( q^{\alpha y}t\right) \left[ x+m%
\right] _{q^{\alpha }}} \\
&=&\left( \sum_{n=0}^{\infty }\left[ y\right] _{q^{\alpha }}^{n}\frac{t^{n}}{%
n!}\right) \left( \sum_{n=0}^{\infty }q^{\alpha \left( n-1\right) y}%
\widetilde{G}_{n,q}^{\left( \alpha ,h\right) }\left( x\right) \frac{t^{n}}{n!%
}\right)
\end{eqnarray*}

by using Cauchy product, we see that%
\begin{equation*}
\sum_{n=0}^{\infty }\left( \sum_{j=0}^{n}\binom{n}{j}q^{\alpha \left(
j-1\right) y}\widetilde{G}_{j,q}^{\left( \alpha ,h\right) }\left( x\right) %
\left[ y\right] _{q^{\alpha }}^{n-j}\right) \frac{t^{n}}{n!}\text{.}
\end{equation*}

Thus, by comparing the coefficients of $\frac{t^{n}}{n!}$, we state the
following corollary.

\begin{corollary}
The following equality holds true:%
\begin{equation}
\widetilde{G}_{n,q}^{\left( \alpha ,h\right) }\left( x+y\right)
=\sum_{j=0}^{n}\binom{n}{j}q^{\alpha \left( j-1\right) y}\widetilde{G}%
_{j,q}^{\left( \alpha ,h\right) }\left( x\right) \left[ y\right] _{q^{\alpha
}}^{n-j}\text{.}  \label{equation 10}
\end{equation}
\end{corollary}

By using Theorem 2.3 and (\ref{equation 10}), we readily derive the
following symmetric relation after some applications.

\begin{theorem}
\label{Theorem 2.5}The following equality holds true:%
\begin{gather*}
\left[ 2\right] _{q^{b}}\sum_{i=0}^{m}\binom{m}{i}\left[ a\right]
_{q^{\alpha }}^{i-1}\left[ b\right] _{q^{\alpha }}^{m-i}\widetilde{G}%
_{i,q^{a}}^{\left( \alpha ,h\right) }\left( bx\right) \widetilde{S}%
_{m-i:q^{b},h+i-1}^{\left( \alpha \right) }\left( a\right) \\
=\left[ 2\right] _{q^{a}}\sum_{i=0}^{m}\binom{m}{i}\left[ b\right]
_{q^{\alpha }}^{i-1}\left[ a\right] _{q^{\alpha }}^{m-i}\widetilde{G}%
_{i,q^{b}}^{\left( \alpha ,h\right) }\left( ax\right) \widetilde{S}%
_{m-i:q^{a},h+i-1}^{\left( \alpha \right) }\left( b\right)
\end{gather*}%
where $\widetilde{S}_{m:q,i}^{\left( \alpha \right) }\left( a\right)
=\sum_{j=0}^{a-1}\left( -1\right) ^{j}q^{ji}\left[ j\right] _{q^{\alpha
}}^{m}$.
\end{theorem}

When $q\rightarrow 1$ into Theorem \ref{Theorem 2.5}, it leads to the
following corollary.

\begin{corollary}
The following identity holds true:%
\begin{gather*}
\sum_{i=0}^{m}\binom{m}{i}a^{i-1}b^{m-i}G_{i}\left( bx\right) S_{m-i}\left(
a\right) \\
=\sum_{i=0}^{m}\binom{m}{i}b^{i-1}a^{m-i}G_{i}\left( ax\right) S_{m-i}\left(
b\right)
\end{gather*}%
where $S_{m}\left( a\right) =\sum_{j=0}^{a-1}\left( -1\right) ^{j}j^{m}$ and 
$G_{n}\left( x\right) $ are called the ordinary Genocchi polynomials which
is defined via the following generating function%
\begin{equation*}
\sum_{n=0}^{\infty }G_{n}\left( x\right) \frac{t^{n}}{n!}=\frac{2t}{e^{t}+1}%
e^{xt}\text{.}
\end{equation*}
\end{corollary}

%%%%%%%%%%%%%%%%%%%%%%%%%%%%%%%%%%%%%%%%%%%%%%%%%%%%%%%%%%%%%%%%

%%%%%%%%%%%%%%%%%%%%%%%%%%%%%%%%%%%%%%%%%%%%%%%%%%%%%%%%%%%%%%%%%%%

\end{document}